\documentclass[12pt]{amsart}
\usepackage{amsfonts}
\usepackage{amsmath}
\usepackage{amssymb}
\usepackage{enumerate}
\usepackage{graphicx}

\makeatletter
\@namedef{subjclassname@2010}{  \textup{2010} Mathematics Subject Classification}
\makeatother

\theoremstyle{definition}

\numberwithin{equation}{section}
\input{tcilatex}

\begin{document}
\title[Note on Poisson cohomology on Weil bundles]{Note on Poisson
cohomology on Weil bundles}
\author[Basile Bossoto]{Basile Guy Richard Bossoto}
\address{Department of Mathematics\\
Faculty of Sciences and Technology, Marien Ngouabi University\\
and \\
Institut de Recherche en Sciences Exactes et Naturelles (IRSEN)\\
Box: 69, Brazzaville, Congo}
\email{bossotob@yahoo.fr}
\author[Mabiala Mikanou]{Olivier Mabiala Mikanou}
\address{Department of Mathematics\\
Faculty of Sciences and Technology, Marien Ngouabi University\\
Box: 69, Brazzaville, Congo }
\email{stive.elg@gmail.com}
\author[Ngu\'{e}ngu\'{e} Louvouandou]{Ap\'{e}p\'{e} Jugend\`{e}ne Ngu\'{e}ngu%
\'{e} Louvouandou}
\address{Department of Mathematics\\
Faculty of Sciences and Technology, Marien Ngouabi University\\
Box: 69, Brazzaville, Congo }
\email{nguengueapepe@gmail.com}
\date{}

\begin{abstract}
The purpose of this note on Poisson cohomology is to show that, if $M$ is a
Poisson manifold and if $A$ is a Weil algebra, then the Weil bundle $M^{A}$
is a Poisson manifold. Therefore, we establish a relationship between
Poisson cohomology $H_{\{.,.\}_{A}}^{\bullet }(M^{A},A)$ with values in $A$
and Poisson cohomology $H_{\{.,.\},\mathbb{R}}^{\bullet }(M^{A})$ with real
values.
\end{abstract}

\subjclass[2010]{Primary 53D17 ; Secondary 58A20.}
\keywords{Weil algebra, Weil bundle, Poisson manifolds, Cohomology.}
\maketitle

\section{Introduction}

In this work, $M$ denotes a smooth manifold of dimension $n$, $C^{\infty
}(M) $ the algebra of smooth functions on $M$ and $A$\ a Weil algebra, i.e a
real commutative, associative, unitary algebra of finite dimension, of the
form $A=%
\mathbb{R}
\oplus \mathfrak{m}$, where $\mathfrak{m}$ is a nilpotent ideal.

We recall that an infinitely near point to $x\in M$\ of type $A$\ is a
homomorphism of $\mathbb{R}$-algebras 
\begin{equation*}
\xi :C^{\infty }(M)\longrightarrow A
\end{equation*}%
such that $\xi (f)-f(x)\in \mathfrak{m}$, for all $f\in C^{\infty }(M)$ \cite%
{wei}. We denote by $M_{x}^{A}$\ the set of all infinitely near points to $%
x\in M$\ of type $A$\ and by 
\begin{equation*}
M^{A}=\bigcup_{x\in M}M_{x}^{A}
\end{equation*}%
the manifold of near points of $M$ of type $A$ or simply the Weil bundle of $%
M$ of type $A$.

We denote by $\pi _{M}:M^{A}\longrightarrow M$, the projection which assigns
every infinitely near point to $x\in M$ to its origin $x$. The triple $%
(M^{A},\pi _{M},M)$ defines the bundle of infinitely near points or simply
Weil bundle \cite{wei}, \cite{mor}, \cite{kms}.

We denote by $C^{\infty }(M^{A})$ the algebra of smooth functions on $M^{A}$
(with values in $%
\mathbb{R}
$) and $C^{\infty }(M^{A},A)$ the algebra of smooth functions on $M^{A}\ $%
with values in $A$. The $C^{\infty }(M)$-module of vector fields on $M$ is
denoted by $\mathfrak{X}(M)$. If $M$ is provided with a Poisson structure
with bracket $\{\cdot ,\cdot \}$, we establish an isomorphism between the
representations from $C^{\infty }(M)$ into $C^{\infty }(M^{A},A)$, and from $%
C^{\infty }(M^{A},A)$ into $C^{\infty }(M^{A},A)$, respectively, defined in 
\cite{nko1}.

When both $M$ and $N$ are smooth manifolds and when $\ h:M\longrightarrow N$
is a smooth mapping, then the mapping%
\begin{equation}
h^{A}:M^{A}\longrightarrow N^{A},\xi \longmapsto h^{A}(\xi ),  \notag
\end{equation}%
such that, for any $g\in C^{\infty }(N)$, 
\begin{equation}
\left[ h^{A}(\xi )\right] (g)=\xi (g\circ h)  \notag
\end{equation}%
is also smooth. When $h$ is a diff{e}omorphism, it is the same for $h^{A}$.

The set $C^{\infty }(M^{A},A)$ of smooth functions on $M^{A}$ with values in 
$A$ is a commutative, unitary algebra over $A$. By identifying $\mathbb{R}%
^{A}$ with $A$, for $f\in C^{\infty }(M)$, the mapping%
\begin{equation}
f^{A}:M^{A}\longrightarrow A,\xi \longmapsto \xi (f)  \notag
\end{equation}%
is differentiable and the mapping%
\begin{equation}
C^{\infty }(M)\longrightarrow C^{\infty }(M^{A},A),f\longmapsto f^{A}\text{,}
\notag
\end{equation}%
is an injective homomorphism of algebras and we have:%
\begin{equation}
(f+g)^{A}=f^{A}+g^{A};(\lambda \cdot f)^{A}=\lambda \cdot f^{A};(f\cdot
g)^{A}=f^{A}\cdot g^{A}  \notag
\end{equation}%
for $\lambda \in \mathbb{R},$ $f,g\in C^{\infty }(M)$.\newline

We denote by $\mathfrak{X}(M^{A})$, the set of all vector fields on $M^{A}.$
According to \cite{bos}, we have the following equivalent assertions:

\begin{enumerate}
\item A vector field on $M^{A}$ is a differentiable section of the tangent
bundle $TM^{A}$.

\item A vector field on $M^{A}$ is a derivation of $C^{\infty }(M^{A})$.

\item A vector field on $M^{A}$ is a linear mapping $X:C^{\infty
}(M)\longrightarrow C^{\infty }(M^{A},A)$ such that 
\begin{equation*}
X(f\cdot g)=X(f)\cdot g^{A}+f^{A}\cdot X(g),\quad \text{for any}\,f,g\in
C^{\infty }(M)\text{.}
\end{equation*}
\end{enumerate}

Thus, the set $\mathfrak{X}(M^{A})$ of vector fields on $M^{A}$ considered
as derivations of $C^{\infty }(M)$ into $C^{\infty }(M^{A},A)$ is a module
over $C^{\infty }(M^{A},A)$.

When%
\begin{equation*}
\theta :C^{\infty }(M)\longrightarrow C^{\infty }(M)
\end{equation*}%
is a vector field on $M$, then the mapping%
\begin{equation*}
\theta ^{A}:C^{\infty }(M)\longrightarrow C^{\infty }(M^{A},A),f\longmapsto 
\left[ \theta (f)\right] ^{A}
\end{equation*}%
is a vector field on $M^{A}$. We say that the vector field $\theta ^{A}$ is
the prolongation to $M^{A}$ of the vector field $\theta $.

If $X$ is a vector field on $M^{A}$, considered as a derivation of $%
C^{\infty }(M)$ into $C^{\infty }(M^{A},A)$, then there exists,\cite{bos}, a
unique derivation%
\begin{equation*}
\widetilde{X}:C^{\infty }(M^{A},A)\longrightarrow C^{\infty }(M^{A},A)
\end{equation*}%
such that:

$1/$ $\widetilde{X}$ is $A$-linear;

$2/$ $\widetilde{X}\left[ C^{\infty}(M^{A})\right] \subset C^{\infty}(M^{A})$%
;

$3/\ \widetilde{X}(f^{A})=X(f)$ for any $f\in C^{\infty}(M)$.

The mapping%
\begin{equation*}
\left[ \cdot ,\cdot \right] :\mathfrak{X}(M^{A})\times \mathfrak{X}%
(M^{A})\longrightarrow \mathfrak{X}(M^{A}),(X,Y)\longmapsto \widetilde{X}%
\circ Y-\widetilde{Y}\circ X\text{,}
\end{equation*}%
is $A$-bilinear and defines a structure of $A$-Lie algebra on $\mathfrak{X}%
(M^{A})$.

Let $(a_{\alpha })_{\alpha =1,...,r}$ be a basis of $A$ and $(a_{\alpha
}^{\ast })_{\alpha =1,...,r}$ be the dual basis.

If%
\begin{equation*}
Y:C^{\infty }(M^{A},A)\longrightarrow C^{\infty }(M^{A},A)
\end{equation*}%
is an $A$-linear derivation such that 
\begin{equation*}
Y(f^{A})=\ \widetilde{X}(f^{A})
\end{equation*}%
for any $f\in C^{\infty }(M)$, then 
\begin{equation}
Y\left[ C^{\infty }(M^{A})\right] \subset C^{\infty }(M^{A})
\end{equation}%
since 
\begin{equation*}
Y(a_{\alpha }^{\ast }\circ f^{A})=\ \widetilde{X}(a_{\alpha }^{\ast }\circ
f^{A})\in C^{\infty }(M^{A})
\end{equation*}%
for any $\alpha =1,2,..,r$. Thus, $Y=\ \widetilde{X}$.

\subsection{Poisson structure on Weil bundle}

A Poisson structure on a smooth manifold $M$ is defined as a bilinear
bracket operation $\left\{ \cdot ,\cdot \right\} $ on $C^{\infty }(M)$ such
that the pair $\left( C^{\infty }(M),\left\{ \cdot ,\cdot \right\} \right) $
is a real Lie algebra and%
\begin{equation*}
\left\{ f,g\cdot h\right\} =\left\{ f,g\right\} \cdot h+g\cdot \left\{
f,h\right\}
\end{equation*}%
for any $f,g,h\in C^{\infty }(M)$. We say in this case that $C^{\infty }(M)$
is a Poisson algebra and $M$ is a Poisson manifold \cite{vai},\cite{lic},%
\cite{car}.

For any $f\in C^{\infty }(M)$, the mapping%
\begin{equation*}
\mathrm{ad}(f):C^{\infty }(M)\longrightarrow C^{\infty }(M),g\longmapsto
\left\{ f,g\right\} ,
\end{equation*}%
is a vector field on $M$ and, for any $g\in C^{\infty }(M)$, we have%
\begin{equation*}
\mathrm{ad}(fg)=f\cdot \mathrm{ad}(g)+g\cdot \mathrm{ad}(f)\text{.}
\end{equation*}

For any $f\in C^{\infty }(M)$, let%
\begin{equation*}
\left[ \mathrm{ad}(f)\right] ^{A}:C^{\infty }(M)\longrightarrow C^{\infty
}(M^{A},A),g\longmapsto \left\{ f,g\right\} ^{A},
\end{equation*}%
be the prolongation of the vector field $\mathrm{ad}(f)$ and let%
\begin{equation*}
\widetilde{\left[ \mathrm{ad}(f)\right] ^{A}}:C^{\infty
}(M^{A},A)\longrightarrow C^{\infty }(M^{A},A)
\end{equation*}%
be the unique $A$-linear derivation such that%
\begin{align*}
\widetilde{\left[ \mathrm{ad}(f)\right] ^{A}}(g^{A})& =\left[ \mathrm{ad}(f)%
\right] ^{A}(g) \\
& =\left\{ f,g\right\} ^{A}
\end{align*}%
for any $g\in C^{\infty }(M)$.

For any $\varphi \in C^{\infty }(M^{A},A)$, the mapping%
\begin{equation*}
\tau _{\varphi }:C^{\infty }(M)\longrightarrow C^{\infty
}(M^{A},A),f\longmapsto -\widetilde{[ad(f)]^{A}}(\varphi ),
\end{equation*}%
is a vector field on $M^{A}$. We denote by%
\begin{equation*}
\widetilde{\tau _{\varphi }}:C^{\infty }(M^{A},A)\longrightarrow C^{\infty
}(M^{A},A)
\end{equation*}%
the unique $A$-linear derivation such that%
\begin{equation*}
\widetilde{\tau _{\varphi }}(f^{A})=\tau _{\varphi }(f)
\end{equation*}%
for any $f\in C^{\infty }(M)$. We have%
\begin{equation*}
\widetilde{\tau }_{f^{A}}=\widetilde{\left[ ad(f)\right] ^{A}}
\end{equation*}%
and for all $\varphi ,\psi \in C^{\infty }(M^{A},A)$ and for $a\in A$: 
\begin{equation*}
\widetilde{\tau _{\varphi +\psi }}=\widetilde{\tau _{\varphi }}+\widetilde{%
\tau _{\psi }}\text{; }\widetilde{\tau _{a\cdot \varphi }}=a\cdot \widetilde{%
\tau _{\varphi }}\text{; }\widetilde{\tau _{\varphi \cdot \psi }}=\varphi
\cdot \widetilde{\tau _{\psi }}+\psi \cdot \widetilde{\tau _{\varphi }}\text{%
.}
\end{equation*}

For any $\varphi ,\psi \in C^{\infty }(M^{A},A)$, if we let%
\begin{equation*}
\left\{ \varphi ,\psi \right\} _{A}=\widetilde{\tau _{\varphi }}(\psi )\text{%
,}
\end{equation*}%
this bracket $\left\{ \cdot ,\cdot \right\} _{A}$ defines a structure of
Poisson $A$-algebra on $C^{\infty }(M^{A},A)$ i.e the mapping 
\begin{equation*}
\left\{ \cdot ,\cdot \right\} _{A}:C^{\infty }(M^{A},A)\times C^{\infty
}(M^{A},A)\longrightarrow C^{\infty }(M^{A},A),(\varphi ,\psi )\longmapsto
\left\{ \varphi ,\psi \right\} _{A}\text{,}
\end{equation*}%
is $A$-bilinear and skew symmetric and for any $\varphi ,\psi ,\phi \in
C^{\infty }(M^{A},A)$,

\begin{enumerate}
\item 
\begin{equation*}
\left\{ \varphi ,\psi \cdot \phi \right\} _{A}=\left\{ \varphi ,\psi
\right\} _{A}\cdot \phi +\psi \cdot \left\{ \varphi ,\phi \right\} _{A}\text{%
.}
\end{equation*}

\item 
\begin{equation*}
\left\{ \varphi ,\left\{ \psi ,\phi \right\} _{A}\right\} _{A}+\left\{ \psi
,\left\{ \phi ,\varphi \right\} _{A}\right\} _{A}+\left\{ \phi ,\left\{
\varphi ,\psi \right\} _{A}\right\} _{A}=0\text{.}
\end{equation*}
\end{enumerate}

We say that the structure of $A$-Poisson manifold on $M^{A}$\ defined by $%
\left\{ \cdot ,\cdot \right\} _{A}$ is the prolongation on $M^{A}$\ of the
structure of Poisson manifold on $M$ defined by $\left\{ \cdot ,\cdot
\right\} $\cite{bos1}.

\section{Poisson cohomology on Weil bundle}

In \cite{oka}, E. Okassa showed that when $\varphi \in A^{\ast }$ is a $%
\mathbb{R}
$-linear form on $A$ and when the manifold $M$ is equipped with a symplectic
structure $(M,\omega )$ ( a pseudo-Riemannian structure $(M,g)$,
respectively), the lift $\varphi \circ \omega ^{A}$ (respectively $\varphi
\circ g^{A}$) on $M^{A}$ is a symplectic structure (a pseudo-Riemannian
structure, repectively) if and only if the dimension of the annihilator of
the maximal ideal $\mathfrak{m}$, ann($\mathfrak{m})$, is equal to $1$ and $%
\varphi \left( \text{ann}(\mathfrak{m})\right) \neq 0$. In \cite{shu}, V. V.
Shurygin studied the Lifts of Poisson structures to $M^{A}$ when $A$ is a
Frobenius algebra. In this paper, we show that if $M$ is a Poisson manifold,
then $M^{A}$ is a Poisson manifold.

Let $(M,\{.,.\})$\ be a Poisson manifold,%
\begin{equation*}
\mathrm{ad}~:f\in C^{\infty }(M)\longmapsto \{f,\cdot \}\in \mathrm{Der}_{%
\mathbb{R}
}(C^{\infty }(M))
\end{equation*}%
the adjoint representation of this structure and $d_{\{.,.\}}$%
\begin{equation*}
\cdots \longrightarrow \bigwedge_{\{.,.\},\mathbb{R}}^{\bullet }\overset{%
d_{\{.,.\}}}{\longrightarrow }\bigwedge_{\{.,.\},\mathbb{R}}^{\bullet
+1}\longrightarrow \cdots 
\end{equation*}%
the Lichnerowicz cohomology operator associated to the representation $%
\mathrm{ad}$ from the Lie-Poisson algebra $(C^{\infty }(M),\{.,.\})$ into $(%
\mathrm{Der}_{%
\mathbb{R}
}(C^{\infty }(M)),[.,.])$ \cite{vai}, where $\bigwedge_{\{.,.\},\mathbb{R}%
}^{\bullet }=$ $\bigwedge_{\{.,.\},\mathbb{R}}^{p}(C^{\infty }(M),C^{\infty
}(M))$ denotes the set of $p$-multilinear skew symmetric mappings from $%
\underset{p-times}{\underbrace{C^{\infty }(M)\times ...\times C^{\infty }(M)}%
}\ $into $C^{\infty }(M)$.

It is well known that if the algebra $C^{\infty }(M)$ is a Poisson algebra
with bracket $\{.,.\}$ (which is $%
\mathbb{R}
$-bilinear), then the algebra $C^{\infty }(M^{A},A)$ is a Poisson algebra
with bracket $\{.,.\}_{A}$ (which is $A$-bilinear) \cite{bos1}. The main
object of this work is to define a $%
\mathbb{R}
$-bilinear Poisson bracket on $C^{\infty }(M^{A})=C^{\infty }(M^{A},%
\mathbb{R}
)$ i.e we show that if $(M,\{.,.\})$ is a Poisson manifold then $M^{A}$ is a
Poisson manifold.

If $\mathrm{Der}_{A}(C^{\infty }(M^{A},A))$ denotes the set of derivations
of $C^{\infty }(M^{A},A)$ which are $A$-linear, the mapping 
\begin{equation}
\tau ~:C^{\infty }(M)\longrightarrow \mathrm{Der}_{A}(C^{\infty
}(M^{A},A)),f\longmapsto \tau (f)  \label{representation}
\end{equation}%
such that 
\begin{equation*}
\tau (f)(\varphi )=\tau _{\varphi }(f)\text{ for all }\varphi \in C^{\infty
}(M^{A},A)
\end{equation*}%
is a $\mathbb{R}$-representation from $(C^{\infty }(M),\{.,.\})$ into $%
(C^{\infty }(M^{A},A),\{.,.\}_{A}~)$. \vskip1mm \noindent This definition
induces (step by step) the construction of the cohomology operator $%
\boldsymbol{d}$ associated with this representation, defined as follows: for
all $p=0,1,...$, for any element $\omega ^{\ast }\in \bigwedge_{\{.,.\}_{A},%
\mathbb{R}}^{p}(C^{\infty }(M),C^{\infty }(M^{A},A))$, for all functions $%
f_{1},...,f_{p},f_{p+1}\in C^{\infty }(M)$, we have 
\begin{equation}
\begin{split}
& \boldsymbol{d}\omega ^{\ast
}(f_{1},...,f_{p+1})=\sum\limits_{j=1}^{p+1}(-1)^{j-1}\tau (f_{j})[\omega
^{\ast }(f_{1},...,\widehat{f_{j}},...,f_{p+1})] \\
& +\sum\limits_{1\leq k<\ell \leq p+1}(-1)^{k+\ell }\omega ^{\ast
}(\{f_{k},f_{\ell }\},f_{1},...,\widehat{f_{k}},...,\widehat{f_{\ell }}%
,...,f_{p+1})
\end{split}%
\end{equation}%
(the notation $\widehat{\cdot }$ means the omission of the component thus
marked). This defines an element of $\bigwedge_{\{.,.\}_{A},\mathbb{R}%
}^{p+1}(C^{\infty }(M),C^{\infty }(M^{A},A))$.

We denote by 
\begin{equation*}
H_{\{.,.\},\mathbb{R}}^{p}(M^{A},A),\qquad p=0,...,n
\end{equation*}%
the associated $\mathbb{R}$-cohomology algebras. As $\dim M=n$, these $%
\mathbb{R}$-vector spaces are zero as soon as $p\geq n+1$.

We can also consider the mapping%
\begin{equation}
\tilde{\tau}~:\varphi \in C^{\infty }(M^{A},A)\longmapsto \widetilde{\tau
_{\varphi }}\in \mathrm{Der}_{A}(C^{\infty }(M^{A},A))
\label{representationbis}
\end{equation}%
as an $A$-representation from the Lie algebra $(C^{\infty
}(M^{A},A),\{.,.\}_{A})$ into the Lie algebra $(C^{\infty
}(M^{A},A),\{.,.\}_{A})$ and in the same way as before, we associate a
cohomology operator%
\begin{equation*}
\boldsymbol{d}_{A}~:\bigwedge_{\{.,.\}_{A},A}^{\bullet }(C^{\infty
}(M^{A},A),C^{\infty }(M^{A},A))\longrightarrow
\bigwedge_{\{.,.\}_{A},A}^{\bullet +1}(C^{\infty }(M^{A},A),C^{\infty
}(M^{A},A))\text{.}
\end{equation*}%
We define for all $\boldsymbol{\omega }^{\ast }\in
\bigwedge_{\{.,.\}_{A},A}^{p}\mathrm{Hom}_{A}(C^{\infty }(M,A),C^{\infty
}(M,A))$, for all functions $\varphi _{1},...,\varphi _{p+1}\in C^{\infty
}(M^{A},A)$: 
\begin{equation}
\begin{split}
& \boldsymbol{d}\boldsymbol{\omega }^{\ast }(\varphi _{1},...,\varphi
_{p+1})=\sum\limits_{j=1}^{p+1}(-1)^{j-1}\widetilde{\tau _{\varphi }}[%
\boldsymbol{\omega }^{\ast }(\varphi _{1},...,\widehat{\varphi _{j}}%
,...,\varphi _{p+1})] \\
& +\sum\limits_{1\leq k<\ell \leq p+1}(-1)^{k+\ell }\boldsymbol{\omega }%
^{\ast }(\{\varphi _{k},\varphi _{\ell }\}_{A},\varphi _{1},...,\widehat{%
\varphi _{k}},...,\widehat{\varphi _{\ell }},...,\varphi _{p+1}) \\
& =\sum\limits_{j=1}^{p+1}(-1)^{j-1}\big\{{\varphi _{j}},\boldsymbol{\omega }%
^{\ast }(\varphi _{1},...,\widehat{\varphi _{j}},...,\varphi _{p+1})\}_{A} \\
& +\sum\limits_{1\leq k<\ell \leq p+1}(-1)^{k+\ell }\boldsymbol{\omega }%
^{\ast }(\{\varphi _{k},\varphi _{\ell }\}_{A},\varphi _{1},...,\widehat{%
\varphi _{k}},...,\widehat{\varphi _{\ell }},...,\varphi _{p+1})\text{.}
\end{split}%
\end{equation}%
The bracket $\{.,.\}_{A}$ above is of course understood as extended in the
sense of the Schouten-Nijenhuis bracket. The associated cohomology $A$%
-algebras will be denoted this time 
\begin{equation*}
H_{\{.,.\}_{A},A}^{p}(M^{A},A)\,,\qquad p=0,...,n\text{.}
\end{equation*}%
(Indeed $M^{A}$ is considered here with its $A$-manifold structure, these
cohomology algebras are null as soon as $p>n+1$).

\begin{theorem}
For all $p\in N$, the cohomology vector space over $%
\mathbb{R}
$, $H_{\{.,.\},\mathbb{R}}^{p}(M^{A},A)$ also has a structure of $A$-module
and coincides with the cohomolgy module over $A$, $H_{\{.,.%
\}_{A},A}^{p}(M^{A},A)$.
\end{theorem}

\begin{proof}
Consider the representation $\tau $ from $(C^{\infty }(M),\{.,.\})$ into $%
(C^{\infty }(M^{A},A),\{.,.\}_{A})$ defined by \eqref{representation} and 
\begin{equation*}
\rho ~:(C^{\infty }(M^{A},A),\{\cdot ,\cdot \}_{A})\rightarrow (\mathrm{Der}%
_{A}(C^{\infty }(M^{A},A)),[.,.])
\end{equation*}%
a representation of Poisson algebra $(C^{\infty }(M^{A},A),\{.,.\}_{A})$. If
we assume that 
\begin{equation*}
\rho (f^{A})=\tau (f)\quad \forall \,f\in C^{\infty }(M),
\end{equation*}%
we observe then that the representations $\rho $ and$\ \tilde{\tau}$
(defined by \eqref{representationbis} coincide for all $\varphi =f^{A}$, $%
f\in C^{\infty }(M)$, therefore are equal since the functions $f^{A}$
generate $C^{\infty }(M^{A},A)$ as $A$-module and that $\rho $ and $\tilde{%
\tau}$ are homomorphisms of $A$-modules. There is therefore a canonical
isomorphism linking the representations $\tau $ and $\tilde{\tau}$ and which
made it possible to construct the two complexes whose cohomology we want to
compare: to $\tau $, we associate $\tilde{\tau}$ as we saw in the previous
section; given on the other hand a representation: $\rho ~:C^{\infty
}(M^{A},A)\rightarrow Der_{A}(C^{\infty }(M^{A},A))$, we associate\textbf{\ }%
the representation%
\begin{equation*}
\underline{\rho }~:C^{\infty }(M)\rightarrow Der_{A}(C^{\infty
}(M^{A},A)),f\ \longmapsto \rho (f^{A}),
\end{equation*}%
thus defining a representation of $(C^{\infty }(M),\{.,.\})$\ in $(C^{\infty
}(M^{A},A),[.,.])$. If $2$\ elements $\omega _{0}^{\ast }$\ and $\omega
_{1}^{\ast }$\ of $\bigwedge_{\{.,.\}_{A},A}^{p}(C^{\infty
}(M^{A},A),C^{\infty }(M^{A},A))$\ are cohomologous relative to the operator 
$d_{A}$, we find by testing them on $(f_{1}^{A},...,f_{p}^{A})$\ that the $2$%
-elements 
\begin{equation*}
\lbrack (f_{1},...,f_{p})\longmapsto \boldsymbol{\omega }_{j}^{\ast
}(f_{1}^{A},...,f_{p}^{A})\big]\in \bigwedge_{\{.,.\}_{A},\mathbb{R}%
}^{p}(C^{\infty }(M),C^{\infty }(M^{A},A)),\ j=0,1,
\end{equation*}%
are cohomologous with respect to the cohomology operator $d$.

Conversely, given two elements $\alpha _{0}^{\ast },$ $\alpha _{1}^{\ast }$ $%
\in \bigwedge_{\{.,.\}_{A},\mathbb{R}}^{p}(C^{\infty }(M),C^{\infty
}(M^{A},A))$ assumed to be cohomologous with respect to the operator $%
\boldsymbol{d}$. The elements $\omega _{0}^{\ast }$ and $\omega _{1}^{\ast }$
defined by extending by $A$-linearity 
\begin{equation*}
(f_{1}^{A},...,f_{p}^{A})\longmapsto \alpha _{j}^{\ast
}(f_{1},...,f_{p})\quad (f_{1},...,f_{p}\in C^{\infty }(M)),\text{ }j=0,1,
\end{equation*}%
are elements of $\bigwedge_{\{.,.\}_{A},A}^{\bullet }(C^{\infty
}(M^{A},A),C^{\infty }(M^{A},A))$ which are cohomologous via this time the
operator $\boldsymbol{d}_{A}$.
\end{proof}

\noindent\ The cohomology defined so will be denoted by the following $%
H_{\{.,.\}_{A}}^{\bullet }(M^{A},A)$.

\begin{proposition}
The restriction of the bracket $\{.,.\}_{A}$ to $C^{\infty }(M^{A})\times
C^{\infty }(M^{A})$, is a Poisson bracket.
\end{proposition}

\begin{proof}
For all $\varphi \in C^{\infty }(M^{A},A)$, the vector field $\tau _{\varphi
}$ is $A$-linear by construction. Thus $\widetilde{\tau _{\varphi }}\left[
C^{\infty }(M^{A})\right] \subset C^{\infty }(M^{A})$. If $F$ and $G$ are in 
$C^{\infty }(M^{A})$, then $\widetilde{\tau _{F}}(G)=\{F,G\}_{A}\in
C^{\infty }(M^{A})$. We define thus a Poisson bracket on $C^{\infty }(M^{A})$
by 
\begin{equation*}
C^{\infty }(M^{A})\times C^{\infty }(M^{A})\longrightarrow C^{\infty
}(M^{A}),(F,G)\longmapsto \{F,G\}_{A}=\widetilde{\tau _{F}}(G)\text{.}
\end{equation*}
\end{proof}

We have the following result:

\begin{theorem}
Let $(M,\{.,.\})$ be a Poisson structure on $M$ and $A$ a Weil algebra. Then 
$M^{A}$ is a Poisson manifold. In addition, we have%
\begin{equation*}
H_{\{.,.\}_{A}}^{\bullet }(M^{A},A)=A\otimes _{%
\mathbb{R}
}H_{\{.,.\},\mathbb{R}}^{\bullet }(M^{A},%
\mathbb{R}
).
\end{equation*}
\end{theorem}

\begin{proof}
According to the previous Proposition, the restriction of the bracket $%
\{.,.\}_{A}$ to $C^{\infty }(M^{A})\times C^{\infty }(M^{A})$ is a Poisson
bracket. Thus, $M^{A}$ is a Poisson manifold with bracket $\{.,.\}_{A}$
restricted to $C^{\infty }(M^{A})\times C^{\infty }(M^{A})$.

When $\{e_{1},...,e_{k}\}$ denotes a basis of $A$ as $\mathbb{R}$-module of
dimension $k$, we now introduce canonical $A$-isomorphism%
\begin{equation*}
\sigma ~:\varphi \in C^{\infty }(M^{A},A)\longmapsto \sigma (\varphi )\in
A\otimes _{%
\mathbb{R}
}C^{\infty }(M^{A})
\end{equation*}%
associating with $\varphi ~:x\mapsto \underset{j=1}{\overset{k}{\dsum }}%
\varphi _{j}(x)\,e_{j}$ the element $\varphi _{j}\otimes e_{j}$, where $%
\varphi =\underset{j=1}{\overset{k}{\dsum }}\varphi _{j}\,e_{j}$ and $%
\varphi _{j}\in C^{\infty }(M^{A},%
\mathbb{R}
)$.

Any $p$-form $\boldsymbol{\omega }^{\ast }\in
\bigwedge_{\{.,.\}_{A}}^{p}(C^{\infty }(M^{A},A),C^{\infty }(M^{A},A))$ of
degree $p$ is equivalent via this isomorphism to the given of $%
\bigwedge^{p}(A\otimes _{%
\mathbb{R}
}C^{\infty }(M^{A}),A\otimes _{%
\mathbb{R}
}C^{\infty }(M^{A}))$. By passing to the quotient relative to the cohomology
induced by the operator $\boldsymbol{d}_{A}$, we deduce the isomorphisms
between the desired cohomology groups. For this purpose, we associate with
the class of $\boldsymbol{\eta }^{\ast }$ the class of%
\begin{equation*}
(\tilde{f}_{1},...,\tilde{f}_{p})\in (A\otimes C^{\infty
}(M^{A}))^{p}\longmapsto \sigma \circ \boldsymbol{\eta }^{\ast }\big(\sigma
^{-1}(\tilde{f}_{1}),...,\sigma ^{-1}(\tilde{f}_{p})\big).
\end{equation*}%
This completes the Proof.
\end{proof}

\subsection*{Acknowledgements}

The authors thank Professor Eug\`{e}ne Okassa and Professor Alain Yger for
their comments, remarks and suggestions as well as International
Mathematical Union (IMU) for the scholarship granted to Olivier Mabiala
Mikanou.




{\normalsize \baselineskip=17pt }

{\normalsize 
}

\end{document}